\documentclass[11pt]{article}
\usepackage{graphicx}
\usepackage{amsmath}
\usepackage{amssymb}
\usepackage{amsthm}
\usepackage[titletoc]{appendix}
\usepackage{url}
\newcommand\al\alpha
\newcommand\be\beta
\newcommand\de\delta
\newcommand\ep\varepsilon
\newcommand\tha\theta
\newcommand\ka\kappa
\newcommand\la\lambda
\newcommand\om\omega

\newcommand\iy\infty
\newcommand\pa\partial

\newcommand{\hyp}[5]{\,\mbox{}_{#1}F_{#2}\!\left(\genfrac{}{}{0pt}{}{#3}{#4};#5\right)}

\renewcommand\Re{\operatorname{Re}}
\numberwithin{equation}{section}

\newtheorem{theorem}{Theorem}

\newtheorem{Remark}[theorem]{Remark}

\begin{document}
	
\title{A remarkable double integral of the product of two Gaussian hypergeometric functions.}
\author{Enno Diekema \footnote{email address: e.diekema@gmail.com}}
\maketitle

\begin{abstract}
\noindent
In this paper a double integral containing two Gaussian hypergeometric functions is discussed. The integral is not found in the literature and a direct computation is not (yet) possible. Therefore, a complete different integral is computed by two methods. The first method gives a combination of hypergeometric functions and the requested original integral, while with the second method this different integral can be computed directly. Equalization of the results from both methods gives the requested integral.
\end{abstract}

\section{Introduction}
\setlength{\parindent}{0cm}
\noindent
During the research of the author about the Rosenblatt distribution function efforts have been made to find a general formula for the cumulants of this function. Computation of the first five cumulants are successfully. These can all be written as combinations of hypergeometric functions. A general formula has not been found (yet). All of this will be described in a subsequent paper. During the research an integral appears that was not known in the literature. It turned out to be very difficult. The integral is
\[
\int_0^1\int_0^1 x^{3-3d}y^{1-d}(1-x\, y)^{-d}\hyp21{1,d}{2-d}{x}\hyp21{1,d}{2-d}{y}dxdy
\]
with $0 \leq d <\, 1$. It is easy to show that this double integral is convergent except for $d=0.8$ because for this value the integral does not exist. 

A direct computation (for example with series expansion) was not successfully. That is the reason the computation has been done in a seemingly cumbersome way. It has started with the following integral which had its origin in \cite{1}
\begin{equation}
\int_0^1\int_0^{y_2}(1-y_3)^{-d}\int_0^{y_3}(y_2-y_4)^{-d}(1-y_4)^{-d}
\int_0^{y_4}(y_2-y_5)^{-d}(y_3-y_5)^{-d}dy_5\, dy_4\, dy_3\, dy_2
\label{2.0a}
\end{equation}
with $0 \leq d <\, 1$.
This integral is computed by two methods. The first method gives a relation to the requested integral. Comparing the two results then give the desired result.

In Section 2 the main result is given. In Section 3 the first method to compute the integral \eqref{2.0a} is treated. In Section 4 the second method is treated. In section 5 the proof is finished. In the Appendix there are a lot of formulas used in this paper including some Thomae transformations.
\vspace{0.3 cm}

\textbf{Remark 1:}\  In the whole paper integrals and summations are interchanged when needed. In all cases it can be shown that this is allowed by the dominant convergence theorem, but the proofs here are omitted while they are irrelevant for the course of the content of the paper.
\vspace{0.3 cm}

\textbf{Remark 2:}\  The simple integrals can usually be found in \cite{2}. Many integrals can also be determined with Mathematica.
\vspace{0.3 cm}
 
\textbf{Remark 3:}\  The properties of the used Pochhammer symbols and the Gamma functions are listed in the Appendix of \cite{4}.
\vspace{0.3 cm}

\textbf{Remark 4:}\  An overview of some Thomae transformations is given in the Appendix of \cite{5}.

\section{Main result}
The main result of this paper is
\begin{align}
\int_0^1\int_0^1 &x^{3-3d}y^{1-d}(1-x\, y)^{-d}\hyp21{1,d}{2-d}{x}\hyp21{1,d}{2-d}{y}dxdy= \nonumber \\[1ex]
&=\dfrac{\Gamma(1-d)}{2(4-5d)}\dfrac{\Gamma(1-d)}{\Gamma(2-2d)}
\hyp32{2-2d,1,2d-1}{3-2d,2-d}{1}- \nonumber\\[1ex]
&-\dfrac{\Gamma(1-d)}{2(4-5d)}
\left(\dfrac{\Gamma(3-3d)}{\Gamma(4-4d)}-\dfrac{\Gamma(1-d)}{\Gamma(2-2d)}\right)
\hyp32{2-2d,1,d}{3-2d,2-d}{1}- \nonumber\\[1ex]
&-\dfrac{\Gamma(1-d)}{2(4-5d)}\dfrac{\Gamma(3-3d)}{\Gamma(4-4d)}
\hyp43{1,d,2-2d,3-3d}{2-d,3-2d,4-4d}{1}
\label{2.00}
\end{align}

\

The integral is convergent for $0 \leq d <1$ except for $d=0.8$ because the integral does not exist for this value.  May be someone can simplify the main result or find another formula.

\section{First method}
Starting with \eqref{2.0a}
\[
\int_0^1\int_0^{y_2}(1-y_3)^{-d}\int_0^{y_3}(y_2-y_4)^{-d}(1-y_4)^{-d}
\int_0^{y_4}(y_2-y_5)^{-d}(y_3-y_5)^{-d}dy_5\, dy_4\, dy_3\, dy_2
\]
and doing the substitutions $y_5=y_4\, z_4,\ y_4=y_3\, z_3,\ y_3=y_2\, z_2,\, y_2=z_2$ gives
\begin{multline}
A=\int_0^1\int_0^1\int_0^1\int_0^1(z_2)^{3-3d}(z_3)^{2-d}
(1-z_2\, z_3)^{-d}(z_4)(1-z_3\, z_4)^{-d}\\
(1-z_2\, z_3\, z_4)^{-d}(1-z_3\, z_4\, z_5)^{-d}(1-z_4\, z_5)^{-d}dz_5\, dz_4\, dz_3\, dz_2
\label{2.1}
\end{multline}
Using \eqref{A.1} for the integral of $z_5$ gives
\begin{align*}
\int_0^1(1-z_3\, z_4\, z_5)^{-d}(1-(z_4\, z_5)^{-d}dz_5
&=\dfrac{1}{(d-1)}(z_4)^{-1}(1-z_4)^{1-d}(1-z_3\, z_4)^{-d}\hyp21{1,d}{2-d}{\dfrac{z_3(1-z_4)}{1-z_3\, z_4}} \\
&-\dfrac{1}{(d-1)}(z_4)^{-1}\hyp21{1,d}{2-d}{z_3}
\end{align*}
Substitution in \eqref{2.1} gives
\begin{align*}
A&=\dfrac{1}{(d-1)}\int_0^1\int_0^1\int_0^1(z_2)^{3-3d}(z_3)^{2-d}
(1-z_2\, z_3)^{-d}(1-z_3\, z_4)^{-2d}
(1-z_2\, z_3\, z_4)^{-d}(1-z_4)^{1-d} \\
&\qquad\qquad\qquad\qquad\qquad\qquad\qquad\qquad\qquad\qquad\qquad\qquad \hyp21{1,d}{2-d}{\dfrac{z_3(1-z_4)}{1-z_3\, z_4}} dz_4\, dz_3\, dz_2 \\
&-\dfrac{1}{(d-1)}\int_0^1\int_0^1\int_0^1(z_2)^{3-3d}(z_3)^{2-d}
(1-z_2\, z_3)^{-d}(1-z_3\, z_4)^{-d}(1-z_2\, z_3\, z_4)^{-d} \\
&\qquad\qquad\qquad\qquad\qquad\qquad\qquad\qquad\qquad\qquad\qquad\qquad\qquad\qquad
\hyp21{1,d}{2-d}{z_3}  dz_4\, dz_3\, dz_2
\end{align*}
After writing the hypergeometric function in the first term as a summation and interchanging the summation and the integrals the result is
\begin{align}
&I_1=\dfrac{1}{(d-1)}\sum_{k=0}^\infty\dfrac{(d)_k}{(2-d)_k}
\int_0^1\int_0^1\int_0^1(z_2)^{3-3d}(z_3)^{2-d+k}(1-z_2\, z_3)^{-d}
(1-z_3\, z_4)^{-2d-k} \nonumber \\
&\qquad\qquad\qquad\qquad\qquad\qquad\qquad\qquad\qquad\qquad
(1-z_2\, z_3\, z_4)^{-d}(1-z_4)^{1-d+k} dz_4\, dz_3\, dz_2 \\
&I_2=\dfrac{1}{(d-1)}\int_0^1\int_0^1\int_0^1(z_2)^{3-3d}(z_3)^{2-d}
(1-z_2\, z_3)^{-d}(1-z_3\, z_4)^{-d}(1-z_2\, z_3\, z_4)^{-d} \nonumber \\
&\qquad\qquad\qquad\qquad\qquad\qquad\qquad\qquad\qquad\qquad\qquad\qquad\qquad\qquad
\hyp21{1,d}{2-d}{z_3}  dz_4\, dz_3\, dz_2
\label{3.3}
\end{align}
where $(.)_k$ is the usual Pochhammer symbol and $A=I_1-I_2$. 

\

At first $I_1$ is treated. The integral of the variable $z_4$ gives with \eqref{A.2a}
\[
\int_0^1 (1-z_3\, z_4)^{-2d-k}(1-z_2\, z_3\, z_4)^{-d} (1-z_4)^{1-d+k} dz_4=
\dfrac{1}{(2-d+k)}
F_1\left(
\begin{array}{c}
1;2d+k,d \\
3-d+k
\end{array}
;z_3,z_2\, z_3\right)
\]
The $F_1$ function is the first Appell function \cite[Sections 5.7-5.14]{3}. Using $\Gamma(a+k)=\Gamma(a)(a)_k$ for
\[
\dfrac{1}{(2-d+k)}=\dfrac{\Gamma(2-d+k)}{\Gamma(3-d+k)}=
\dfrac{\Gamma(2-d)}{\Gamma(3-d)}\dfrac{(2-d)_k}{(3-d)_k}
\]
writing the Appell function as a double summation \eqref{A.2b} and interchanging the summations and the integrals gives
\begin{multline}
I_1=\dfrac{1}{(d-1)}\dfrac{\Gamma(2-d)}{\Gamma(3-d)}\sum_{k=0}^\infty
\dfrac{(d)_k}{(3-d)_k}\sum_{i=0}^\infty\sum_{j=0}^\infty
\dfrac{(1)_{i+j}(2d+k)_i(d)_j}{(3-d+k)_{i+j}}\dfrac{1}{i!j!} \\
\int_0^1\int_0^1(z_2)^{3-3d+j}(z_3)^{2-d+i+j+k}(1-z_2\, z_3)^{-d}dz_2\, dz_3
\label{I1a}
\end{multline}

The double integral can be computed with standard methods.
\begin{multline*}
\int_0^1\int_0^1(z_2)^{3-3d+j}(z_3)^{2-d+i+j+k}(1-z_2\, z_3)^{-d}dz_2\, dz_3= \\
\dfrac{\Gamma(1-d)}{(2d-1+i+k)}\left(\dfrac{\Gamma(4-3d+j)}{\Gamma(5-4d+j)}-
\dfrac{\Gamma(3-d+i+j+k)}{\Gamma(4-2d+i+j+k)}\right)
\end{multline*}
Application to \eqref{I1a} and rearranging gives
\begin{align*}
I_1&=\dfrac{\Gamma(1-d)}{(d-1)}\dfrac{\Gamma(2-d)}{\Gamma(3-d)}
\dfrac{\Gamma(4-3d)}{\Gamma(5-4d)}\dfrac{\Gamma(2d-1)}{\Gamma(2d)}
\sum_{k=0}^\infty\dfrac{(d)_k}{(3-d)_k}\dfrac{(2d-1)_k}{(2d)_k}
\sum_{j=0}^\infty\dfrac{(d)_j}{(3-d+k)_j}\dfrac{(4-3d)_j}{(5-4d)_j} \\
&\qquad\qquad\qquad\qquad\qquad\qquad\qquad\qquad\qquad\qquad\qquad\qquad\qquad\qquad
\sum_{i=0}^\infty\dfrac{(1+j)_i(2d-1+k)_i}{(3-d+j+k)_i}\dfrac{1}{i!} \\
&-\dfrac{\Gamma(1-d)}{(d-1)}\dfrac{\Gamma(2-d)}{\Gamma(4-2d)}
\dfrac{\Gamma(2d-1)}{\Gamma(2d)}\sum_{k=0}^\infty\dfrac{(d)_k(2d-1)_k}{(4-2d)_k(2d)_k}
\sum_{j=0}^\infty\dfrac{(d)_j}{(4-2d+k)_j}
\sum_{i=0}^\infty\dfrac{(1+j)_i(2d-1+k)_i}{(4-2d+j+k)_i}\dfrac{1}{i!}
\end{align*}
The last summations in both terms can be written as $_2F_1$ hypergeometric functions with unit argument and are both convergent. Using
\begin{equation}
\hyp21{a,b}{c}{1}=\dfrac{\Gamma(c)\Gamma(c-a-b)}{\Gamma(c-a)\Gamma(c-b)}
\label{3.00}
\end{equation}
with $c-a-b>0$ gives
\begin{align*}
I_1&=\dfrac{\Gamma(1-d)}{(d-1)}\dfrac{\Gamma(2-d)}{\Gamma(3-d)}
\dfrac{\Gamma(4-3d)}{\Gamma(5-4d)}\dfrac{\Gamma(2d-1)}{\Gamma(2d)}
\sum_{k=0}^\infty\dfrac{(d)_k}{(3-d)_k}\dfrac{(2d-1)_k}{(2d)_k}
\sum_{j=0}^\infty\dfrac{(d)_j}{(3-d+k)_j}\dfrac{(4-3d)_j}{(5-4d)_j} \\
&\qquad\qquad\qquad\qquad\qquad\qquad\qquad\qquad\qquad\qquad\qquad\qquad\qquad\qquad
\dfrac{\Gamma(3-d+j+k)\Gamma(3-3d)}{\Gamma(2-d+k)\Gamma(4-3d+j)} \\
&-\dfrac{\Gamma(1-d)}{(d-1)}\dfrac{\Gamma(2-d)}{\Gamma(4-2d)}
\dfrac{\Gamma(2d-1)}{\Gamma(2d)}\sum_{k=0}^\infty\dfrac{(d)_k(2d-1)_k}{(4-2d)_k(2d)_k}
\sum_{j=0}^\infty\dfrac{(d)_j}{(4-2d+k)_j}
\dfrac{\Gamma(4-2d+j+k)\Gamma(4-4d)}{\Gamma(3-2d+k)\Gamma(5-4d+j)}
\end{align*}
Simplification gives
\begin{align*}
I_1&=\dfrac{\Gamma(1-d)}{(d-1)}\dfrac{\Gamma(3-3d)}{\Gamma(5-4d)}
\dfrac{\Gamma(2d-1)}{\Gamma(2d)}
\sum_{k=0}^\infty\dfrac{(d)_k(2d-1)_k}{(2-d)_k(2d)_k}
\sum_{j=0}^\infty\dfrac{(1)_j(d)_j}{(5-4d)_j}\dfrac{1}{j!} \\
&-\dfrac{\Gamma(1-d)}{(d-1)}\dfrac{\Gamma(2-d)}{\Gamma(3-2d)}
\dfrac{\Gamma(2d-1)}{\Gamma(2d)}\dfrac{\Gamma(4-4d)}{\Gamma(5-4d)}
\sum_{k=0}^\infty\dfrac{(d)_k(2d-1)_k}{(3-2d)_k(2d)_k}
\sum_{j=0}^\infty\dfrac{(1)_j(d)_j}{(5-4d)_j}\dfrac{1}{j!}
\end{align*}
The last summation is well known. Application of \eqref{3.00} gives
\begin{align*}
I_1&=\dfrac{\Gamma(1-d)}{(d-1)}\dfrac{\Gamma(2d-1)}{\Gamma(2d)}
\dfrac{\Gamma(3-3d)\Gamma(4-5d)}{\Gamma(4-4d)\Gamma(5-5d)}
\sum_{k=0}^\infty\dfrac{(d)_k(2d-1)_k}{(2-d)_k(2d)_k} \\
&-\dfrac{\Gamma(1-d)}{(d-1)}\dfrac{\Gamma(2-d)}{\Gamma(3-2d)}
\dfrac{\Gamma(2d-1)}{\Gamma(2d)}\dfrac{\Gamma(4-5d)}{\Gamma(5-5d)}
\sum_{k=0}^\infty\dfrac{(d)_k(2d-1)_k}{(3-2d)_k(2d)_k}
\end{align*}
The summations can be written as hypergeometric functions. Then at last there is
\begin{align*}
I_1&=\dfrac{\Gamma(1-d)}{(d-1)}\dfrac{\Gamma(2d-1)}{\Gamma(2d)}
\dfrac{\Gamma(3-3d)\Gamma(4-5d)}{\Gamma(4-4d)\Gamma(5-5d)}
\hyp32{1,d,2d-1}{2-d,2d}{1}\nonumber \\
&+\dfrac{\Gamma(1-d)^2}{\Gamma(3-2d)}\dfrac{\Gamma(2d-1)}{\Gamma(2d)}
\dfrac{\Gamma(4-5d)}{\Gamma(5-5d)}\hyp32{1,d,2d-1}{3-2d,2d}{1}
\end{align*}
For reasons that become apparent in the last sector the transformation \eqref{A.7} will be used for both hypergeometric functions. The result is
\begin{align}
I_1&=\dfrac{\Gamma(1-d)}{(1-d)(2-2d)}\dfrac{\Gamma(3-3d)\Gamma(4-5d)}{\Gamma(4-4d)\Gamma(5-5d)}\hyp32{2-2d,1,d}{3-2d,2-d}{1}- \nonumber \\
&-\dfrac{\Gamma(1-d)^3}{\Gamma(3-2d)\Gamma(2-d)}\dfrac{\Gamma(4-5d)}{\Gamma(5-5d)}
\hyp32{2-2d,1,2d-1}{3-2d,2-d}{1}
\label{I1}
\end{align}

\

\

Next $I_2$ is treated. Integration to the variable $z_4$ in \eqref{3.3} using \eqref{A.1} gives
\begin{align*}
&\int_0^1(1-z_3\, z_4)^{-d}(1-z_2\, z_3\, z_4)^{-d}dz_4= \\
&=\dfrac{1}{(d-1)}(z_3)^{-1}(1-z_3)^{1-d}(1-z_2\, z_3)^{-d}
\hyp21{1,d}{2-d}{\dfrac{z_2(1-z_3)}{1-z_z\, z_3}}-
\dfrac{1}{(d-1)}(z_3)^{-1}\hyp21{1,d}{2-d}{z_2}
\end{align*}
Substitution in \eqref{3.3} gives
\begin{align}
I_2&=\dfrac{1}{(d-1)^2}\int_0^1\int_0^1(z_2)^{3-3d}(z_3)^{1-d}(1-z_3)^{1-d}(1-z_2\, z_3)^{-2d}
\hyp21{1,d}{2-d}{\dfrac{z_2(1-z_3)}{1-z_z\, z_3}} \nonumber \\
&\qquad\qquad\qquad\qquad\qquad\qquad\qquad\qquad\qquad\qquad\qquad\qquad\qquad
\hyp21{1,d}{2-d}{z_3} dz_3\, dz_2- \label{3.5} \\
&-\dfrac{1}{(d-1)^2}\int_0^1\int_0^1(z_2)^{3-3d}(z_3)^{1-d}(1-z_2\, z_3)^{-d} 
\hyp21{1,d}{2-d}{z_2}\hyp21{1,d}{2-d}{z_3} dz_3\, dz_2 \nonumber
\end{align}
The required integral appears in the second term. Applying \eqref{A.2} to the first hypergeometric function in the first term gives after some simplification
\begin{align*}
\hyp21{1,d}{2-d}{\dfrac{z_2(1-z_3)}{1-z_z\, z_3}}
&=\dfrac{\Gamma(2-d)\Gamma(2d-1)}{\Gamma(d)}(z_2)^{d-1}(1-z_2)^{1-2d}(1-z_3)^{d-1}(1-z_2\, z_3)^d+ \\
&+\dfrac{\Gamma(1-2d)\Gamma(2-d)}{\Gamma(2-2d)\Gamma(1-d)}\hyp21{1,d}{2d}
{1-\dfrac{z_2(1-z_3)}{1-z_2\, z_3}}
\end{align*}
Application of the transformation \eqref{A.0} gives
\begin{align*}
\hyp21{1,d}{2-d}{\dfrac{z_2(1-z_3)}{1-z_z\, z_3}}
&=\dfrac{\Gamma(2-d)\Gamma(2d-1)}{\Gamma(d)}(z_2)^{d-1}(1-z_2)^{1-2d}(1-z_3)^{d-1}(1-z_2\, z_3)^d+ \\
&+\dfrac{\Gamma(1-2d)\Gamma(2-d)}{\Gamma(2-2d)\Gamma(1-d)}
\left(\dfrac{1-z_2}{1-z_2\, z_3}\right)^{d-1}
\hyp21{d,2d-1}{2d}{\dfrac{1-z_2}{1-z_2\, z_3}}
\end{align*}
Substitution in \eqref{3.5} gives after some simplification
\begin{align}
&I_2(a)=\dfrac{1}{(d-1)^2}\dfrac{\Gamma(2-d)\Gamma(2d-1)}{\Gamma(d)}
\int_0^1\int_0^1(z_2)^{2-2d}(1-z_2)^{1-2d}(z_3)^{1-d}(1-z_2\, z_3)^{-d} \nonumber \\
&\qquad\qquad\qquad\qquad\qquad\qquad\qquad\qquad\qquad\qquad\qquad\qquad\qquad\qquad
\hyp21{1,d}{2-d}{z_3}dz_3\, dz_2 \label{3.6} \\
&I_2(b)=\dfrac{1}{(d-1)}\dfrac{\Gamma(1-2d)}{\Gamma(2-2d)}\int_0^1\int_0^1(z_2)^{2-2d}
(z_3)^{1-d}(1-z_2\, z_3)^{1-3d}\hyp21{d,2d-1}{2d}{\dfrac{1-z_2}{1-z_2\, z_3}} \nonumber \\
&\qquad\qquad\qquad\qquad\qquad\qquad\qquad\qquad\qquad\qquad\qquad\qquad\qquad
\qquad\hyp21{1,d}{2-d}{z_3}dz_3\, dz_2 \label{3.7} \\
&I_2(c)=\dfrac{1}{(d-1)^2}\int_0^1\int_0^1(z_2)^{3-3d}(z_3)^{1-d}(1-z_2\, z_3)^{-d} 
\hyp21{1,d}{2-d}{z_2}\hyp21{1,d}{2-d}{z_3} dz_3\, dz_2 \nonumber
\end{align}
with $I_2=I_2(a)-I_2(b)-I_2(c)$.

\

At first the computation of $I_2(a)$ is treated. From the double integral the integration to $z_2$ gives.
\[
\int_0^1(z_2)^{2-2d}(1-z_2)^{1-2d}(1-z_2\, z_3)^{-d} dz_2=
\dfrac{\Gamma(3-2d)\Gamma(2-2d)}{\Gamma(5-4d)}\hyp21{d,3-2d}{5-4d}{z_3}
\]
Substitution in the double integral gives
\begin{multline*}
A=\int_0^1\int_0^1(z_2)^{2-2d}(1-z_2)^{1-2d}(z_3)^{1-d}(1-z_2\, z_3)^{-d}
\hyp21{1,d}{2-d}{z_3}dz_3\, dz_2= \\
=\dfrac{\Gamma(3-2d)\Gamma(2-2d)}{\Gamma(5-4d)}\int_0^1(z_3)^{1-d}
\hyp21{1,d}{2-d}{z_3}\hyp21{d,3-2d}{5-4d}{z_3}dz_3
\end{multline*}
Application of the transformation \eqref{A.0} on the first hypergeometric function and changing the integration variable $z_3$ into $1-z_3$ gives
\begin{multline*}
A=\dfrac{\Gamma(3-2d)\Gamma(2-2d)}{\Gamma(5-4d)}\int_0^1(z_3)^{1-2d}(1-z_3)^{1-d}
\hyp21{1-d,2-2d}{2-d}{1-z_3} \\
\hyp21{d,3-2d}{5-4d}{1-z_3}dz_3
\end{multline*}
Application of \eqref{A.3} gives
\begin{align*}
A&=\dfrac{\Gamma(3-2d)\Gamma(2-2d)}{\Gamma(5-4d)}
\dfrac{\Gamma(5-4d)\Gamma(2-3d)}{\Gamma(5-5d)\Gamma(3-2d)}
\hyp32{1,d,2-2d}{2-d,3d-1}{1}+ \\
&+\dfrac{\Gamma(3-2d)\Gamma(2-2d)}{\Gamma(5-4d)}
\dfrac{\Gamma(4-5d)\Gamma(5-4d)\Gamma(3-3d)\Gamma(2-d)\Gamma(3d-2)}
{\Gamma(5-5d)\Gamma(4-4d)\Gamma(3-2d)\Gamma(d)}\hyp21{4-5d,2-2d}{4-4d}{1}
\end{align*}
Because the $3F2$ hypergeometric function is not convergent for $0 \leq d\leq 1$ a Thomae transformation \eqref{A.4} is applied. The $2F1$ hypergeometric function can be evaluated using \eqref{3.00}. After lot of manipulations the result is
\begin{multline*}
A=\dfrac{\Gamma(2-2d)\Gamma(2-3d)}{\Gamma(5-5d)}
\dfrac{\Gamma(2-d)\Gamma(3d-2)}{\Gamma(d)^2}\hyp32{2-2d,2d-1,3d-2}{d,3d-1}{1}+ \\
+\dfrac{\Gamma(4-5d)\Gamma(3-3d)\Gamma(2-d)\Gamma(3d-2)^2}{\Gamma(5-5d)\Gamma(d)^2}
\end{multline*}
Substitution in \eqref{3.6} gives for  $I_2(a)$
\begin{multline*}
I_2(a)=\dfrac{\Gamma(2d-1)\Gamma(1-d)^2\Gamma(3d-2)}{\Gamma(d)^3}
\dfrac{\Gamma(2-2d)\Gamma(2-3d)}{\Gamma(5-5d)}\hyp32{2-2d,2d-1,3d-2}{d,3d-1}{1}+ \\
\dfrac{\Gamma(2d-1)\Gamma(1-d)^2\Gamma(4-5d)\Gamma(3-3d)\Gamma(3d-2)^2}{\Gamma(5-5d)\Gamma(d)^3}
\end{multline*}
Application of \eqref{A.7} gives again after lots of manipulation
\begin{equation}
I_2(a)=\dfrac{\Gamma(1-d)^2\Gamma(2-2d)^2}{\Gamma(d)^2}
\dfrac{\Gamma(2d-1)\Gamma(3-3d)}{\Gamma(5-5d)(4-5d)\Gamma(4-4d)}
\hyp32{2-2d,3-3d,4-5d}{4-4d,5-5d}{1}
\label{I2a}
\end{equation}

\

Next the computation of $I_2(b)$. In \eqref{3.7} a modification of the first hypergeometric function is done setting
\[
z_1=\dfrac{1-z_2}{1-z_2\, z_3}
\]Then the integral of $I_2(b)$ becomes
\[
B=\int_0^1\int_0^1(1-z_1)^{2-2d}(z_2)^{1-d}(1-z_2)^{2-3d}(1-z_1\, z_2)^{5d-5}
\hyp21{d,2d-1}{2d}{z_1}\hyp21{1,d}{2-d}{z_2}dz_1\, dz_2
\]
The hypergeometric functions can be written as summations. Interchanging the summations and the integrals gives
\begin{multline*}
B=\sum_{k=0}^\infty\dfrac{(d)_k(2d-1)_k}{(2d)_k}\dfrac{1}{k!}
\sum_{j=0}^\infty\dfrac{(d)_j}{(2-d)_j} \\
\int_0^1\int_0^1(z_1)^k(1-z_1)^{2-2d}(z_2)^{1-d+j}(1-z_2)^{2-3d}(1-z_1\, z_2)^{5d-5}
dz_1\, dz_2
\end{multline*}
The double integral can computed by standard methods.
\begin{multline*}
\int_0^1\int_0^1(z_1)^k(1-z_1)^{2-2d}(z_2)^{1-d+j}(1-z_2)^{2-3d}(1-z_1\, z_2)^{5d-5}
dz_1\, dz_2= \\
=\dfrac{\Gamma(1+k)\Gamma(3-3d)\Gamma(3-2d)\Gamma(2-d+j)}{\Gamma(5-4d+j)\Gamma(4-2d+k)}
\hyp32{5-5d,2-d+j,1+k}{5-4d+j,4-2d+k}{1}
\end{multline*}
Substitution in $B$ gives after some simplification
\[
B=\dfrac{\Gamma(3-3d)\Gamma(3-2d)\Gamma(2-d)}{\Gamma(5-4d)\Gamma(4-2d)}
\sum_{k=0}^\infty\dfrac{(d_k)(2d-1)_k}{(2d)_k(4-2d)_k}
\sum_{j=0}^\infty\dfrac{(d)_j}{(5-4d)_j}\hyp32{5-5d,2-d+j,1+k}{5-4d+j,4-2d+k}{1}
\]
To execute the summations the hypergeometric function can be transformed. Using the Thomae transformation \eqref{A.5} gives
\[
B=\dfrac{\Gamma(3-3d)\Gamma(2-d)}{\Gamma(5-4d)}
\sum_{k=0}^\infty\dfrac{(d_k)(2d-1)_k}{(2d)_k(2)_k}
\sum_{j=0}^\infty\dfrac{(d)_j}{(5-4d)_j}\hyp32{3-3d,d+j,1+k}{5-4d+j,2+k}{1}
\]
Writing the hypergeometric function as a summation and interchanging the summations gives
\[
B=\dfrac{\Gamma(3-3d)\Gamma(2-d)}{\Gamma(5-4d)}
\sum_{m=0}^\infty\dfrac{(3-3d)_m}{m!}
\sum_{k=0}^\infty\dfrac{(d)_k(2d-1)_k(1+k)_m}{(2d)_k(2)_k(2+k)_m}
\sum_{j=0}^\infty\dfrac{(d)_j(d+j)_m}{(5-4d)_j(5-4d+j)_m}
\]
The last summations are known. With $(a+i)_j=(a+j)_i(a)_j/(a)_i$ the result is
\[
\sum_{j=0}^\infty\dfrac{(d)_j(d+j)_m}{(5-4d)_j(5-4d+j)_m}=
\dfrac{\Gamma(4-5d)\Gamma(5-4d)}{\Gamma(5-5d)\Gamma(4-4d)}\dfrac{(d)_m}{(2-d)_m} 
\]
and using \eqref{A.6} gives
\[
\sum_{k=0}^\infty\dfrac{(d)_k(2d-1)_k(1+k)_m}{(2d)_k(2)_k(2+k)_m}=
\dfrac{\Gamma(1-d)\Gamma(2-2d)}{\Gamma(3-2d)}\dfrac{\Gamma(2d)}{\Gamma(d)}
\dfrac{(2-2d)_m}{(3-2d)_m}+
\dfrac{(1-2d)}{(2-2d)}\dfrac{\Gamma(1-d)}{\Gamma(2-d)}
\dfrac{(2-2d)_m(1)_m}{(3-2d)_m(2-d)_m}
\]
Application while writing the summations as hypergeometric functions results in
\begin{align*}
B&=\dfrac{\Gamma(1-d)^2}{2}\dfrac{\Gamma(2d)}{\Gamma(d)}
\dfrac{\Gamma(3-3d)\Gamma(4-5d)}{\Gamma(5-5d)\Gamma(4-4d)}
\hyp32{d,2-2d,3-3d}{3-2d,4-4d}{1}+ \\
&+\dfrac{\Gamma(1-d)(1-2d)}{(2-2d)}
\dfrac{\Gamma(3-3d)\Gamma(4-5d)}{\Gamma(5-5d)\Gamma(4-4d)}
\hyp43{1,d,2-2d,3-3d}{2-d,3-2d,4-4d}{1}
\end{align*}
Substitution in the formula for $I_2(b)$ gives at last
\begin{align*}
I_2(b)&=\dfrac{\Gamma(1-d)}{(d-1)}\dfrac{\Gamma(2-2d)}{\Gamma(3-2d)}
\dfrac{\Gamma(3-3d)\Gamma(4-5d)}{\Gamma(5-5d)\Gamma(4-4d)}
\hyp43{1,d,2-2d,3-3d}{2-d,3-2d,4-4d}{1}- \nonumber \\
&-\dfrac{\Gamma(1-d)^2\Gamma(1-2d)}{\Gamma(3-2d)}
\dfrac{\Gamma(2d)}{\Gamma(d)}\dfrac{\Gamma(3-3d)\Gamma(4-5d)}{\Gamma(5-5d)\Gamma(4-4d)}
\hyp32{d,2-2d,3-3d}{3-2d,4-4d}{1}
\end{align*}
Application of the Thomae transformation \eqref{A.51} for the hypergeometric function in the second term results in
\begin{align}
I_2(b)&=\dfrac{\Gamma(1-d)}{(d-1)}\dfrac{\Gamma(2-2d)}{\Gamma(3-2d)}
\dfrac{\Gamma(3-3d)\Gamma(4-5d)}{\Gamma(5-5d)\Gamma(4-4d)}
\hyp43{1,d,2-2d,3-3d}{2-d,3-2d,4-4d}{1}+ \nonumber \\
&+\dfrac{\Gamma(1-d)^2\Gamma(2-2d)^2}{\Gamma(d)^2}
\dfrac{\Gamma(2d-1)\Gamma(3-3d)}{\Gamma(5-5d)(4-5d)\Gamma(4-4d)}
\hyp32{2-2d,3-3d,4-5d}{4-4d,5-5d}{1}
\label{I2b}
\end{align}

\section{Second method}
Repeating \eqref{2.0a} gives
\begin{equation}
\int_0^1\int_0^{y_2}(1-y_3)^{-d}\int_0^{y_3}(y_2-y_4)^{-d}(1-y_4)^{-d}
\int_0^{y_4}(y_2-y_5)^{-d}(y_3-y_5)^{-d}dy_5\, dy_4\, dy_3\, dy_2
\label{2.0}
\end{equation}
Starting with the integral of $dy_5$ and using the following property
\[
\int_0^a\int_0^x f(x,y)dy\, dx=\int_0^a\int_y^a f(x,y)dx\, dy
\]
the integral can be rewritten as
\[
A=\int_0^1\int_0^{y_2}(1-y_3)^{-d}\int_0^{y_3}(y_2-y_5)^{-d}(y_3-y_5)^{-d}
\int_{y_5}^{y_3}(y_2-y_4)^{-d}(1-y_4)^{-d}dy_4\, dy_5\, dy_3\, dy_2
\]
For the integral of $y_4$ the transformation $y_4=(y_3-y_5)z+y_5$ should be done. This gives
\begin{multline*}
\int_{y_5}^{y_3}(y_2-y_4)^{-d}(1-y_4)^{-d}dy_4= \\
=(1-y_5)^{-d}(y_2-y_5)^{-d}(y_3-y_5)
\int_0^1\left(1-\dfrac{y_3-y_5}{y_2-y_5}z\right)^{-d}
\left(1-\dfrac{y_3-y_5}{1-y_5}z\right)^{-d}dz
\end{multline*}
Using \eqref{A.1} gives after a lot of manipulations $A=J_1-J_2$ with
\begin{multline*}
J_1=\dfrac{1}{(1-d)}\int_0^1\int_0^{y_2}\int_0^{y_3}(1-y_3)^{-d}(1-y_5)^{-d}
(y_2-y_5)^{1-2d}(y_3-y_5)^{-d}  \\
\hyp21{1,d}{2-d}{\dfrac{y_2-y_5}{1-y_5}}dy_5\, dy_3\, dy_2
\end{multline*}
\begin{multline*}
J_2=\dfrac{1}{(1-d)}\int_0^1\int_0^{y_2}\int_0^{y_3}(1-y_3)^{-d}(y_2-y_3)^{1-d}
(y_2-y_5)^{-d}(y_3-y_5)^{-d} \\
\hyp21{1,d}{2-d}{\dfrac{y_2-y_3}{1-y_3}}dy_5\, dy_3\, dy_2
\end{multline*}
The integral with the variable $y_5$ in $J_2$ gives
\[
\int_0^{y_3}(y_3-y_5)^{-d}(y_2-y_5)^{-d}dy_5=\dfrac{1}{(1-d)}
\hyp21{1,d}{2-d}{\dfrac{y_3}{y_2}}
\]
Substitution and applying the transformation \eqref{A.0} to the other hypergeometric functions gives
\begin{multline*}
J_1=\dfrac{1}{(1-d)}\int_0^1\int_0^{y_2}\int_0^{y_3}(1-y_3)^{-d}(1-y_5)^{-d}
(y_2-y_5)^{1-2d}(y_3-y_5)^{-d} \\
\hyp21{1,d}{2-d}{\dfrac{y_2-y_5}{1-y_5}}dy_5\, dy_3\, dy_2
\end{multline*}
\begin{multline*}
J_2=\dfrac{1}{(1-d)^2}\int_0^1\int_0^{y_2}(y_2)^{-d}(y_3)^{1-d}(1-y_3)^{-2d}
(y_2-y_3)^{1-d}\hyp21{1,d}{2-d}{\dfrac{y_3}{y_2}} \\
\hyp21{1,d}{2-d}{\dfrac{y_2-y_3}{1-y_3}}dy_3\, dy_2
\end{multline*}

\

The substitutions $y_5=y_4\, z_4,\ y_4=y_3\, z_3,\ y_3=y_2\, z_2,\, y_2=z_2$ give
\begin{multline*}
J_1=\dfrac{1}{(1-d)}\int_0^1\int_0^1\int_0^1(z_2)^{3-3d}(z_3)^{1-d}(1-z_2\, z_3)^{-d}(1-z_3\, z_5)^{1-2d}(1-z_5)^{-d}(1-z_2\, z_3\, z_5)^{-d} \\
\hyp21{1,d}{2-d}{\dfrac{1-z_3\, z_5}{1-z_2\, z_3\, z_5}z_2}dz_5\, dz_3\, dz_2
\end{multline*}
\begin{multline*}
J_2=\dfrac{1}{(1-d)^2}\int_0^1\int_0^1
(z_2)^{3-3d}(z_3)^{1-d}(1-z_2\, z_3)^{-2d}(1-z_3)^{1-d}
\hyp21{1,d}{2-d}{z_3} \\
\hyp21{1,d}{2-d}{\dfrac{1-z_3}{1-z_2\, z_3}z_2}dz_3\, dz_2
\end{multline*}
Starting with $J_1$ and writing the hypergeometric function as a summation gives \begin{multline*}
J_1=\dfrac{1}{(1-d)}\sum_{k=0}^\infty\dfrac{(d)_k}{(2-d)_k}\int_0^1\int_0^1\int_0^1
(z_2)^{3-3d+k}(z_3)^{1-d}(1-z_2\, z_3)^{-d}(1-z_3\, z_5)^{1-2d+k} \\
(1-z_5)^{-d}(1-z_2\, z_3\, z_5)^{-d-k} dz_5\, dz_3\, dz_2
\end{multline*}
Integration to the variable $z_5$ gives with \eqref{A.2a}
\[
\int_0^1(1-z_5)^{-d}(1-z_3\, z_5)^{1-2d+k}(1-z_2\, z_3\, z_5)^{-d-k} dz_5=\dfrac{1}{(1-d)}
F_1\left(
\begin{array}{c}
	1;2d-1-k,d+k \\
	2-d\qquad\qquad
\end{array}
;z_3,z_2\, z_3\right)
\]
Substitution gives
\begin{multline*}
J_1=\dfrac{1}{(1-d)^2}\sum_{k=0}^\infty\dfrac{(d)_k}{(2-d)_k}\int_0^1\int_0^1
(z_2)^{3-3d+k}(z_3)^{1-d}(1-z_2\, z_3)^{-d} \\
F_1\left(
\begin{array}{c}
	1;2d-1-k,d+k \\
	2-d\qquad\qquad
\end{array}
;z_3,z_2\, z_3\right)
dz_3\, dz_2
\end{multline*}
Writing the Appell function  as a double summation \eqref{A.2b} gives
\begin{multline*}
J_1=\dfrac{1}{(1-d)^2}\sum_{k=0}^\infty\sum_{i=0}^\infty \sum_{j=0}^\infty
\dfrac{(d)_k}{(2-d)_k}\dfrac{(1)_{i+j}(2d-1-k)_i(d+k)_j}{(2-d)_{i+j}}\dfrac{1}{i!j!} \\
\int_0^1\int_0^1(z_2)^{3-3d+k+j}(z_3)^{1-d+i+j}(1-z_2\, z_3)^{-d} dz_3\, dz_2
\end{multline*}
The integral can be computed. The result is
\begin{multline*}
\int_0^1\int_0^1(z_2)^{3-3d+k+j}(z_3)^{1-d+i+j}(1-z_2\, z_3)^{-d} dz_3\, dz_2= \\
=\Gamma(1-d)\dfrac{\Gamma(2-2d-i+k)}{\Gamma(3-2d-i+k)}
\dfrac{\Gamma(2-d+i+j)}{\Gamma(3-2d+i+j)}-\Gamma(1-d)
\dfrac{\Gamma(2-2d-i+k)}{\Gamma(3-2d-i+k)}\dfrac{\Gamma(4-3d+j+k)}{\Gamma(5-4d+j+k)}
\end{multline*}
Substitution gives
\begin{align*}
&J_1(a)=\dfrac{\Gamma(1-d)}{(1-d)^2}\sum_{k=0}^\infty\sum_{i=0}^\infty \sum_{j=0}^\infty
\dfrac{(d)_k}{(2-d)_k}\dfrac{(1)_{i+j}(2d-1-k)_i(d+k)_j}{(2-d)_{i+j}}\dfrac{1}{i!j!} \\
&\qquad\qquad\qquad\qquad\qquad\qquad\qquad\qquad\qquad\qquad\qquad
\dfrac{\Gamma(2-2d-i+k)}{\Gamma(3-2d-i+k)}
\dfrac{\Gamma(2-d+i+j)}{\Gamma(3-2d+i+j)} \\
&J_1(b)=\dfrac{\Gamma(1-d)}{(1-d)^2}\sum_{k=0}^\infty\sum_{i=0}^\infty \sum_{j=0}^\infty
\dfrac{(d)_k}{(2-d)_k}\dfrac{(1)_{i+j}(2d-1-k)_i(d+k)_j}{(2-d)_{i+j}}\dfrac{1}{i!j!} \\
&\qquad\qquad\qquad\qquad\qquad\qquad\qquad\qquad\qquad\qquad\qquad
\dfrac{\Gamma(2-2d-i+k)}{\Gamma(3-2d-i+k)}
\dfrac{\Gamma(4-3d+j+k)}{\Gamma(5-4d+j+k)}
\end{align*}
with $J_1=J_1(a)-J_1(b)$. Using $(a)_{i+j}=(a+j)_i(a)_j$ in $J_1(a)$ gives
\begin{multline*}
J_1(a)=\dfrac{\Gamma(1-d)}{(1-d)^2}\dfrac{\Gamma(2-2d)}{\Gamma(3-2d)}
\dfrac{\Gamma(2-d)}{\Gamma(3-2d)}
\sum_{k=0}^\infty\dfrac{(d)_k(2-2d)_k}{(2-d)_k(3-2d)_k}
\sum_{j=0}^\infty\dfrac{(d+k)_j}{(3-2d)_j} \\
\sum_{i=0}^\infty\dfrac{(1+j)_i(2-2d-k)_i}{(3-2d+j)_i}\dfrac{1}{i!}
\end{multline*}
The last summation results in
\[
\sum_{i=0}^\infty\dfrac{(1+j)_i(2-2d-k)_i}{(3-2d+j)_i}\dfrac{1}{i!}=
\dfrac{\Gamma(3-2d+j)\Gamma(4-4d+k)}{\Gamma(2-2d)\Gamma(5-4d+j+k)}
\]
Applying this result and after a lot of manipulation there is
\[
J_1(a)=\dfrac{\Gamma(1-d)}{(1-d)^2}\dfrac{\Gamma(2-d)}{\Gamma(3-2d)}
\dfrac{\Gamma(4-4d)}{\Gamma(5-4d)}
\sum_{k=0}^\infty\dfrac{(d)_k(2-2d)_k}{(2-d)_k(3-2d)_k}
\sum_{j=0}^\infty\dfrac{(d+k)_j}{(5-4d+k)_j}
\]
The last summation applies
\[
\sum_{j=0}^\infty\dfrac{(d+k)_j}{(5-4d+k)_j}=
\dfrac{\Gamma(5-4d)\Gamma(4-5d)}{\Gamma(4-4d)\Gamma(5-5d)}\dfrac{(5-4d)_k}{(4-4d)_k}
\]
After substitution the final result is
\begin{equation}
J_1(a)=\dfrac{\Gamma(1-d)}{(1-d)^2}\dfrac{\Gamma(2-d)}{\Gamma(3-2d)}
\dfrac{\Gamma(4-5d)}{\Gamma(5-5d)}\hyp32{1,d,2-2d}{2-d,3-2d}{1}
\label{I1aa}
\end{equation}

\

Next $J_1(b)$. Using $(a)_{i+j}=(a+j)_i(a)_j$ gives
\begin{multline*}
J_1(b)=\dfrac{\Gamma(1-d)}{(1-d)^2}\dfrac{\Gamma(2-2d)}{\Gamma(3-2d)}
\sum_{k=0}^\infty\dfrac{(d)_k(2-2d)_k}{(2-d)_k(3-2d)_k}\dfrac{\Gamma(4-3d+k)}{\Gamma(5-4d+k)}
\sum_{j=0}^\infty\dfrac{(d+k)_j(4-3d+k)_j}{(2-d)_j(5-4d+k)_j}  \\
\sum_{i=0}^\infty\dfrac{(1+j)_i(2d-2-k)_i}{(2-d+j)_i}\dfrac{1}{i!}
\end{multline*}
The last summation results in
\[
\sum_{i=0}^\infty\dfrac{(1+j)_i(2d-2-k)_i}{(2-d+j)_i}\dfrac{1}{i!}=
\dfrac{\Gamma(2-d+j)\Gamma(3-3d+k)}{\Gamma(1-d)\Gamma(4-3d+j+k)}
\]
Applying this summation and after a lot of manipulation the result is
\[
J_1(b)=\dfrac{\Gamma(1-d)}{(1-d)}\dfrac{\Gamma(2-2d)}{\Gamma(3-2d)}
\dfrac{\Gamma(3-3d)}{\Gamma(5-4d)}
\sum_{k=0}^\infty\dfrac{(d)_k(2-2d)_k(3-3d)_k}{(2-d)_k(3-2d)_k(5-4d)_k}
\sum_{j=0}^\infty\dfrac{(d+k)_j}{(5-4d+k)_j}
\]
The second summation results in
\[
\sum_{j=0}^\infty\dfrac{(d+k)_j}{(5-4d+k)_j}=
\dfrac{\Gamma(5-4d+k)\Gamma(4-5d)}{\Gamma(4-4d+k(\Gamma(5-5d)}
\]
After substitution there remains at last
\begin{equation}
J_1(b)=\dfrac{\Gamma(1-d)}{2(1-d)^2}\dfrac{\Gamma(3-3d)\Gamma(4-5d)}{\Gamma(4-4d)\Gamma(5-5d)}\hyp43{1,d,2-2d,3-3d}{2-d,3-2d,4-4d}{1}
\label{I1bb}
\end{equation}

\

Next the integral for $J_2$ will be computed.
\begin{multline*}
J_2=\dfrac{1}{(1-d)^2}\int_0^1\int_0^1
(z_2)^{3-3d}(z_3)^{1-d}(1-z_2\, z_3)^{-2d}(1-z_3)^{1-d}
\hyp21{1,d}{2-d}{z_3} \\
\hyp21{1,d}{2-d}{\dfrac{1-z_3}{1-z_2\, z_3}z_2}dz_3\, dz_2
\end{multline*}
Writing the hypergeometric functions as summations and interchanging the summations and the integrals gives
\begin{multline*}
J_2=\dfrac{1}{(1-d)^2}\sum_{k=0}^\infty\dfrac{(d)_k}{(2-d)_k}
\sum_{j=0}^\infty\dfrac{(d)_j}{(2-d)_j}\int_0^1(z_3)^{1-d+k}(1-z_3)^{1-d+j} \\
\int_0^1(z_2)^{3-3d+j}(1-z_2\, z_3)^{-2d-j} dz_2\, dz_3
\end{multline*}
The integral with $z_2$ gives a hypergeometric function
\begin{align*}
\int_0^1(z_2)^{3-3d+j}(1-z_2\, z_3)^{-2d-j} dz_2
&=\dfrac{1}{(4-3d+j)}\hyp21{4-3d+j,2d+j}{5-3d+j}{z_3} \\
&=\dfrac{1}{(4-3d+j)}(1-z_3)^{1-2d-j}\hyp21{1,5-5d}{5-3d+j}{z_3}
\end{align*}
The integral with $z_3$ gives
\begin{multline*}
\dfrac{1}{(4-3d+j)}
\int_0^1(z_3)^{1-d+k}(1-z_3)^{2-3d}\hyp21{1,5-5d}{5-3d+j}{z_3}= \\
=\Gamma(3-3d)\dfrac{\Gamma(4-3d+j)}{\Gamma(5-3d+j)}
\dfrac{\Gamma(2-d+k)}{\Gamma(5-4d+k)}
\hyp32{1,5-5d,2-d+k}{5-3d+j,5-4d+k}{1}
\end{multline*}
Substitution in $J_2$ gives after some simplification
\begin{multline*}
J_2=\dfrac{\Gamma(3-3d)}{(1-d)^2}\dfrac{\Gamma(4-3d)}{\Gamma(5-3d)}
\dfrac{\Gamma(2-d)}{\Gamma(5-4d)}
\sum_{k=0}^\infty\dfrac{(d)_k}{(5-4d)_k}
\sum_{j=0}^\infty\dfrac{(d)_j(4-3d)_j}{(2-d)_j(5-3d)_j} \\
\hyp32{1,5-5d,2-d+k}{5-3d+j,5-4d+k}{1}
\end{multline*}
For the hypergeometric function the Thomae transformation \eqref{A.4} will be used.  Application gives after a lot of manipulations
\[
J_2=\dfrac{\Gamma(3-3d)\Gamma(1-d)^2}{\Gamma(5-4d)\Gamma(3-d)}
\sum_{k=0}^\infty\dfrac{(d)_k}{(5-4d)_k}
\sum_{j=0}^\infty\dfrac{(d)_j}{(3-d)_j} \hyp32{1,3-3d,d+k}{3-d+j,5-4d+k}{1}
\]
Writing the hypergeometric function as a summation and after some rearranging the result is
\[
J_2=\dfrac{\Gamma(3-3d)\Gamma(1-d)^2}{\Gamma(5-4d)\Gamma(3-d)}
\sum_{m=0}^\infty(3-3d)_m
\sum_{k=0}^\infty\dfrac{(d)_k}{(5-4d)_k}\dfrac{(d+k)_m}{(5-4d+k)_m}
\sum_{j=0}^\infty\dfrac{(d)_j}{(3-d)_j}\dfrac{1}{(3-d+j)_m}
\]
Tha last two summations can be worked out. 
\begin{multline*}
J_2=\dfrac{\Gamma(3-3d)\Gamma(1-d)^2}{\Gamma(5-4d)\Gamma(3-d)}
\sum_{m=0}^\infty(3-3d)_m
\dfrac{\Gamma(4-5d)\Gamma(4-5d+m)(d)_m}{\Gamma(5-5d)\Gamma(4-4d+m)(5-4d)_m} \\
\dfrac{\Gamma(2-2d+m)\Gamma(3-d+m)}{\Gamma(3-2d+m)\Gamma(2-d+m)(3-d)_m}
\end{multline*}
The final result is
\begin{equation}
J_2=\dfrac{\Gamma(1-d)}{2(1-d)^2}
\dfrac{\Gamma(3-3d)\Gamma(4-5d)}{\Gamma(4-4d)\Gamma(5-5d)}
\hyp43{1,d,2-2d,3-3d}{2-d,3-2d,4-4d}{1}
\label{I2cc}
\end{equation}
This is equal to $J_1(b)$.

\
\section{Final proof}
In this section the results from the previous sections are gathered.

\

\underline{Results of Section 3}.

In Section 3 the result for the integral is $I_1-I_2$. For $I_1$ \eqref{I1} gives
\begin{align*}
I_1&=\dfrac{\Gamma(1-d)}{(1-d)(2-2d)}\dfrac{\Gamma(3-3d)\Gamma(4-5d)}{\Gamma(4-4d)\Gamma(5-5d)}\hyp32{2-2d,1,d}{3-2d,2-d}{1}- \nonumber \\
&-\dfrac{\Gamma(1-d)^3}{\Gamma(3-2d)\Gamma(2-d)}\dfrac{\Gamma(4-5d)}{\Gamma(5-5d)}
\hyp32{2-2d,1,2d-1}{3-2d,2-d}{1}
\end{align*}
There is $I_2=I_2(a)-I_2(b)-I_2(c)$.
For $I_2(a)$ \eqref{I2a} gives
\begin{equation*}
I_2(a)=\dfrac{\Gamma(1-d)^2\Gamma(2-2d)^2}{\Gamma(d)^2}
\dfrac{\Gamma(2d-1)\Gamma(3-3d)}{\Gamma(5-5d)(4-5d)\Gamma(4-4d)}
\hyp32{2-2d,3-3d,4-5d}{4-4d,5-5d}{1}
\end{equation*}
For $I_2(b)$ \eqref{I2b} gives
\begin{align*}
I_2(b)&=\dfrac{\Gamma(1-d)}{(d-1)}\dfrac{\Gamma(2-2d)}{\Gamma(3-2d)}
\dfrac{\Gamma(3-3d)\Gamma(4-5d)}{\Gamma(5-5d)\Gamma(4-4d)}
\hyp43{1,d,2-2d,3-3d}{2-d,3-2d,4-4d}{1}+  \\
&+\dfrac{\Gamma(1-d)^2\Gamma(2-2d)^2}{\Gamma(d)^2}
\dfrac{\Gamma(2d-1)\Gamma(3-3d)}{\Gamma(5-5d)(4-5d)\Gamma(4-4d)}
\hyp32{2-2d,3-3d,4-5d}{4-4d,5-5d}{1}
\end{align*}
$I_2(c)$ is the requested integral. Setting $x=z_2$ and $y=z_3$ gives
\[
I_2(c)=\dfrac{1}{(d-1)^2}\int_0^1\int_0^1 x^{3-3d}y^{1-d}(1-x\, y)^{-d}\hyp21{1,d}{2-d}{x}\hyp21{1,d}{2-d}{y}dxdy
\]

\

\underline{Results of Section 4}.

In Section 4 the result for the integral is $J_1-J_2$. For $J_1$ there is $J_1=J_1(a)-J_1(b)$. For $J_1(a)$ \eqref{I1aa} gives
\[
J_1(a)=\dfrac{\Gamma(1-d)}{(1-d)^2}\dfrac{\Gamma(2-d)}{\Gamma(3-2d)}
\dfrac{\Gamma(4-5d)}{\Gamma(5-5d)}\hyp32{1,d,2-2d}{2-d,3-2d}{1}
\]
For $J_1(b)$ \eqref{I1bb} gives
\[
J_1(b)=\dfrac{\Gamma(1-d)}{2(1-d)^2}
\dfrac{\Gamma(3-3d)\Gamma(4-5d)}{\Gamma(4-4d)\Gamma(5-5d)}
\hyp43{1,d,2-2d,3-3d}{2-d,3-2d,4-4d}{1}
\]
For $J_2$ \eqref{I2cc} gives
\[
J_2=\dfrac{\Gamma(1-d)}{2(1-d)^2}
\dfrac{\Gamma(3-3d)\Gamma(4-5d)}{\Gamma(4-4d)\Gamma(5-5d)}
\hyp43{1,d,2-2d,3-3d}{2-d,3-2d,4-4d}{1}
\]
Setting the result of Section 3 equal to the result of Section 4 gives
\begin{multline*}
I_1-I_2(a)-I_2(b)+\dfrac{1}{(d-1)^2}\int_0^1\int_0^1 x^{3-3d}y^{1-d}(1-x\, y)^{-d}\hyp21{1,d}{2-d}{x}\hyp21{1,d}{2-d}{y}dxdy= \\
J_1(a)-J_1(b)-J_2
\end{multline*}
After substitutions and some rearrangements the final result is \eqref{2.00} and this ends the proof.

\section{Appendix}
In this Appendix there are a lot of formulas used in this paper.

First the standard Gaussian transformations of the hypergeometric function.
\begin{equation}
\hyp21{a,b}{c}{z}=(1-z)^{c-a-b}\hyp21{c-a,c-b}{c}{z}
=(1-z)^{-a}\hyp21{a,c-b}{c}{\dfrac{z}{z-1}}
\label{A.0}
\end{equation}
A well known integral 
\begin{align}
\int_0^1(1-p\, z)^{-d}(1-q\, z)^{-d}dz
&=\dfrac{1}{(1-d)}p^{-1}\left(\dfrac{p}{p-q}\right)^d\hyp21{1-d,d}{2-d}{\dfrac{q}{q-p}}\nonumber \\
&-\dfrac{1}{(1-d)}p^{-1}(1-p)^{-d}\left(\dfrac{p(1-q)}{p-q}\right)^d\hyp21{1-d,d}{2-d}{\dfrac{q(1-p)}{q-p}}
\label{A.1}
\end{align}
\cite[5.8.2.(5)]{3} gives for the Appell $F_1$ function
\begin{multline}
F_1\left(
\begin{array}{c}
	\alpha;\beta,\beta' \\
	\gamma\qquad
\end{array}
;x\, y\right)
=\dfrac{\Gamma(\gamma)}{\Gamma(\alpha)\Gamma(\beta)}
\int_0^1u^{\alpha-1}(1-u)^{\gamma-\alpha-1}(1-u\, x)^{-\beta}(1-u\, y)^{-\beta'}du \\
\Re(\alpha)>0,\ \Re(\gamma-\alpha)>0
\label{A.2a}
\end{multline}
Writing this function as a double summation gives with \cite[5.7.1.(6)]{3}
\begin{equation}
F_1\left(
\begin{array}{c}
	\alpha;\beta,\beta' \\
	\gamma\qquad
\end{array}
;x\, y\right)=
\sum_{i=0}^\infty\sum_{j=0}^\infty
\dfrac{(\alpha)_{i+j}(\beta)_i(\beta')_j}{(\gamma)_{i+j}i!j!}x^i\, y^j
\label{A.2b}
\end{equation}
From \cite[2.10(1)]{3} there is
\begin{align}
\hyp21{a,b}{c}{z}
&=\dfrac{\Gamma(c)\Gamma(c-a-b)}{\Gamma(c-a)\Gamma(c-b)}
\hyp21{a,b}{a+b-c+1}{1-z}+ \nonumber \\
&+\dfrac{\Gamma(c)\Gamma(a+b-c)}{\Gamma(a)\Gamma(b)}(1-z)^{c-a-b}
\hyp21{c-a,c-b}{c-a-b+1}{1-z}
\label{A.2}
\end{align}
From \cite[2.21.9(12)]{2} there is
\begin{align}
&\int_0^1x^{\alpha-1}(1-x)^{c-1}\hyp21{a,b}{c}{1-x}\hyp21{a',b'}{c'}{1-x}dx= \nonumber \\
&=\dfrac{\Gamma(c)\Gamma(c')\Gamma(c'-a'-b')\Gamma)(\alpha)\Gamma(c-a-b+\alpha)}
{\Gamma(c-a+\alpha)\Gamma(c-b+\alpha)\Gamma(c'-a')\Gamma(c'-b')}
\hyp43{a',b',\alpha,c-a-b+\alpha}{c-a+\alpha,c-b+\alpha,a'+b'-c'+1}{1}+ \nonumber \\
&+\dfrac{\Gamma(c)\Gamma(c')\Gamma(a'+b'-\underline{c'})\Gamma(c'-a'-b'+\alpha)\Gamma(c+c'-a-a'-b-b'+\alpha)}{\Gamma(a')\Gamma(b')\Gamma(c+c'-a-a'-b'+\alpha)\Gamma(c+c'-a'-b-b'+\alpha)} \nonumber \\
&\qquad\qquad\qquad
\hyp43{c'-a',c'-b',c'-a'-b'+\alpha,c+c'-a-a'-b-b'+\alpha}{c'-a'-b'+1,c+c'-a-a'-b'+\alpha,c+c'-a'-b-b'+\alpha}{1}
\label{A.3}
\end{align}
Note the little mistake with the underlined $c'$ which is here corrected and communicated with Yu.A. Brychkov.

\

For the $3F2$ hypergeometric functions with unit argument there are a lot of Thomae transformations. Here are three samples used in this paper.
\begin{align}
&\hyp32{a,b,c}{e,f}{1}=\dfrac{\Gamma(e+f-a-b-c)\Gamma(e)}{\Gamma(e-a)\Gamma(e+f-b-c)}
\hyp32{a,f-c,f-b}{e+f-b-c,f}{1}
\label{A.4} \\
&\hyp32{a,b,c}{e,f}{1}=\dfrac{\Gamma(e+f-a-b-c)\Gamma(f)}{\Gamma(f-c)\Gamma(e+f-b-a)}
\hyp32{a,f-c,f-b}{e+f-b-c,f}{1}
\label{A.5} \\
&\hyp32{a,b,c}{e,f}{1}=\dfrac{\Gamma(e+f-a-b-c)\Gamma(f)\Gamma(e)}
{\Gamma(b)\Gamma(e+f-b-c)\Gamma(e+f-a-b)}  \nonumber \\
&\qquad\qquad\qquad\qquad\qquad\qquad\qquad\qquad\qquad\qquad
\hyp32{e+f-a-b-c,f-b,e-b}{e+f-b-c,e+f-a-b}{1}
\label{A.51}
\end{align}

\

From \cite[7.4.4.(16)]{2} there is
\begin{equation}
\hyp32{a,b,c}{a+1,b+1}{1}=\dfrac{\Gamma(1-c)}{(a-b)}
\left(a\dfrac{\Gamma(b+1)}{\Gamma(b-c+1)}-b\dfrac{\Gamma(a+1)}{\Gamma(a-c+1)}\right)
\label{A.6}
\end{equation}
This last theorem can be proved by a Thomae transformation.

\

From \cite[7.4.4.(11)]{2} there is
\begin{align}
\hyp32{a,b,c}{a+1,d}{1}
&=\dfrac{\Gamma(1-c)\Gamma(d)\Gamma(a+1)\Gamma(b-a)}{\Gamma(b)\Gamma(d-a)\Gamma(a-c+1)}- \nonumber \\
&-\dfrac{a\Gamma(1-c)\Gamma(d)}{\Gamma(b-c+1)\Gamma(d-b)(b-a)}\hyp32{b.b-d+1,b-a}{b-a+1,b-c+1}{1}
\label{A.7}
\end{align}

\end{document}